\documentclass{article}
\newtheorem{lemma}{Lemma}
\newtheorem{prop}{Proposition}
\newtheorem{coro}{Corollary}
\newtheorem{defi}{Definition}
\newtheorem{theor}{Theorem}
\usepackage{amsmath}
\usepackage{amssymb}
\newcommand{\D}{\mathbf{D}}
\newcommand{\N}{\mathbb{N}}
\usepackage[english]{babel}
\usepackage{graphicx} 
\usepackage{hyperref}
\usepackage{float}

\title{An elementary proof of the existence and uniqueness of solutions to an initial value problem}
\author{Luca Tanganelli Castrillón}
\date{February 12, 2024}

\begin{document}

\maketitle
\begin{abstract}
    In this note, we show a classical result on the local existence and uniqueness of a solution to an initial value problem subject to a Lipschitz condition. We use only elementary tools from mathematical analysis, without involving any integration. We proceed by showing that the Cauchy iterates converge on a dense subset of the interval and subsequently proving that the extension of this limit function to the whole interval is a solution to the Cauchy problem.
\end{abstract}
\section{The initial value problem}
The initial value problem consists of finding a function $x(t)$ such that
$$(\text{IVP})\begin{cases}
    x'(t)=f(t,x(t)),\qquad t\geq t_0 \\ x(t_0)=x_0
\end{cases}$$
where $f:D\subset\mathbb{R}^2\to\mathbb{R}$, $t_0,x_0\in \mathbb{R}$.

Under certain conditions, the problem has a solution and it is unique.

Let's assume that there exist $a,M,L>0$ such that $U:=[t_0,t_0+a]\times [x_0-Ma,x_0+Ma]\subset D$ and that for all $\mathbf{u},\mathbf{v}\in U$
$$|f(\mathbf{u})-f(\mathbf{v})|\leq L\left\|\mathbf{u}-\mathbf{v}\right\|_1,\qquad |f(\mathbf{u})|\leq M.$$
That is, $f$ is $L$-Lipschitz-continuous and bounded by $M$ in $U$. The interest in the existence of such $U$ lies in the fact that, as a consequence of the mean value theorem, any solution function of (IVP) defined in $[t_0,t_0+a]$ must take values in $[x_0-Ma,x_0+Ma]$ which will facilitate our work.

Then we are going to prove that there exists a function $x:[t_0,t_0+a]\to \mathbb{R}$ that verifies the IVP. Afterwards, we will see that it is the only function defined in $[t_0,t_0+a]$ that verifies (IVP).

For simplicity in notation, we will assume $t_0=0$. The original problem is recovered without complication through a translation.

\begin{defi}
    Given $m\geq 1$, we define $\{x_0^m, x_{\frac{a}{2^m}}^m,\ldots,x_{\frac{ak}{2^m}}^m,\ldots,x_{\frac{a(2^m-1)}{2^m}}^m,x_a^m\}$ as follows.
    $\begin{cases}
        x_0^m=x_0\\
        x_{\frac{a(k+1)}{2^m}}^m=x_{\frac{ak}{2^m}}^m+\frac{a}{2^m}f\left(\frac{ak}{2^m},x_{\frac{ak}{2^m}}^m\right)
    \end{cases}$
\end{defi}
\begin{defi}
    We will call $\mathbf{D}=\{\frac{ak}{2^m}:0\leq k\leq 2^m, k,m\in\mathbb{N}\}$.
\end{defi}
An important property of $\D$ is that it is dense in $[0,a]$.
\begin{lemma}
    For all $c,d\in\D$, if $m\in\N$ is such that $x_c^m$ and $x_d^m$ are well defined, then
    \begin{itemize}
        \item $|x_c^m-x_0|\leq Mc\leq Ma\implies (c,x_c^m)\in U$;
        \item $|x_c^m-x_d^m|\leq M|c-d|$.
    \end{itemize}
\end{lemma}
\textbf{Proof.} The first is by induction on $0\leq k\leq 2^m$. The second is a consequence of the triangular inequality and the bounding of $f$ by $M$ in $U$.

\begin{prop}
    For all $d\in\D$, the sequence $(x_d^m)_{m=m_d}^\infty$ is convergent, where $m_d$ is the first natural for which the numbers $x_d^m$, $m\geq m_d$, are well defined. Moreover, the convergence is uniform in the variable $d$.
\end{prop}

\textbf{Proof.}
For each $m\geq m_d$, let $e_d^m=|x_d^{m+1}-x_d^m|$.
Then observe that $e_0^m=0$. Also, if $k>0$
$$x_{\frac{ak}{2^m}}^m=x_{\frac{a(k-1)}{2^m}}^m+\frac{a}{2^m}f\left(\frac{a(k-1)}{2^m},x_{\frac{a(k-1)}{2^m}}^m\right)$$
$$x_{\frac{ak}{2^m}}^{m+1}=x_{\frac{a(k-1)}{2^m}}^{m+1}+\frac{a}{2^{m+1}}\left[f\left(\frac{a(k-1)}{2^m},x_{\frac{a(k-1)}{2^m}}^{m+1}\right)+f\left(\frac{a(2k-1)}{2^{m+1}},x_{\frac{a(2k-1)}{2^{m+1}}}^{m+1}\right)\right].$$
Subtracting the equation above from the one below and taking modules, the following inequality is obtained.
\begin{equation*}e_\frac{ak}{2^m}^m\leq e_\frac{a(k-1)}{2^m}^m+\frac{a}{2^{m+1}} 
    \Bigg(\left| f  \left(\frac{a(k-1)}{2^m},x_{\frac{a(k-1)}{2^m}}^{m+1}\right)-f  \left(\frac{a(k-1)}{2^m},x_{\frac{a(k-1)}{2^m}}^m\right)\right| + \end{equation*}\begin{equation} \left|f  \left(\frac{a(2k-1)}{2^{m+1}},x_{\frac{a(2k-1)}{2^{m+1}}}^{m+1}\right) -f 
 \left(\frac{a(k-1)}{2^m},x_{\frac{a(k-1)}{2^m}}^m\right)\right|\Bigg).\end{equation}

Now, \begin{equation}\left\|\left(\frac{a(k-1)}{2^m},x_{\frac{a(k-1)}{2^m}}^{m+1}\right)-\left(\frac{a(k-1)}{2^m},x_{\frac{a(k-1)}{2^m}}^m\right)\right\|_1=|x_{\frac{a(k-1)}{2^m}}^{m+1}-x_{\frac{a(k-1)}{2^m}}^m|=e_\frac{a(k-1)}{2^m}^m\end{equation}
then, applying the Lipschitz property,
$$\left| f  \left(\frac{a(k-1)}{2^m},x_{\frac{a(k-1)}{2^m}}^{m+1}\right)-f  \left(\frac{a(k-1)}{2^m},x_{\frac{a(k-1)}{2^m}}^m\right)\right|\leq Le_\frac{a(k-1)}{2^m}^m.$$
Also,
\begin{align*}\left\|\left(\frac{a(2k-1)}{2^{m+1}},x_{\frac{a(2k-1)}{2^{m+1}}}^{m+1}\right)-\left(\frac{a(k-1)}{2^m},x_{\frac{a(k-1)}{2^m}}^m\right)\right\|_1&=\frac{a}{2^{m+1}}+|x_{\frac{a(2k-1)}{2^{m+1}}}^{m+1}-x_{\frac{a(k-1)}{2^m}}^m|\\
&\leq \frac{a}{2^{m+1}}+|x_{\frac{a(2k-1)}{2^{m+1}}}^{m+1}-x_{\frac{a(k-1)}{2^m}}^{m+1}|+|x_{\frac{a(k-1)}{2^m}}^{m+1}-x_{\frac{a(k-1)}{2^m}}^m|\\
&=\frac{a}{2^{m+1}}+\left|\frac{a}{2^{m+1}}f\left(\frac{a(k-1)}{2^m},x_{\frac{a(k-1)}{2^m}}^{m+1}\right)\right|+e_\frac{a(k-1)}{2^m}^m\\
&\leq \frac{a}{2^{m+1}}+\frac{a}{2^{m+1}}M+e_\frac{a(k-1)}{2^m}^m.
\end{align*}
Applying the Lipschitz property again:
\begin{equation}\left|f  \left(\frac{a(2k-1)}{2^{m+1}},x_{\frac{a(2k-1)}{2^{m+1}}}^{m+1}\right) -f 
 \left(\frac{a(k-1)}{2^m},x_{\frac{a(k-1)}{2^m}}^m\right)\right|\leq L\left(\frac{a(M+1)}{2^{m+1}}+e_\frac{a(k-1)}{2^m}^m\right).\end{equation}
Combining (2) and (3) in (1):
\begin{align*}e_\frac{ak}{2^m}^m\leq e_\frac{a(k-1)}{2^m}^m+\frac{a}{2^{m+1}}\left(Le_\frac{a(k-1)}{2^m}^m+L\left(\frac{a(M+1)}{2^{m+1}}+e_\frac{a(k-1)}{2^m}^m\right)\right)&=e_\frac{a(k-1)}{2^m}^m\left(1+\frac{aL}{2^m}\right)+\frac{a^2L(M+1)}{2^{2m+2}}\\
&=e_\frac{a(k-1)}{2^m}^m\left(1+\alpha\right)+\beta
\end{align*}
where $\alpha=\frac{aL}{2^m}$ and $\beta=\frac{a^2L(M+1)}{2^{2m+2}}$. Now, making use of this inequality, it is observed that
\begin{align*}
    e_\frac{ak}{2^m}^m&\leq e_\frac{a(k-1)}{2^m}^m(1+\alpha)+\beta\\
    &\leq (e_\frac{a(k-2)}{2^m}^m(1+\alpha)+\beta)(1+\alpha)+\beta\\&=e_\frac{a(k-2)}{2^m}^m(1+\alpha)^2+\beta(1+(1+\alpha))\\
    &\leq (e_\frac{a(k-3)}{2^m}^m(1+\alpha)+\beta)(1+\alpha)^2+\beta(1+(1+\alpha))\\&=e_\frac{a(k-3)}{2^m}^m(1+\alpha)^3+\beta(1+(1+\alpha)+(1+\alpha)^2)\\
    &\leq \ldots\\
    &\leq e_\frac{a(k-k)}{2^m}^m(1+\alpha)^k+\beta(1+(1+\alpha)+(1+\alpha)^2+\ldots+(1+\alpha)^k)\\
    &=\beta(1+(1+\alpha)+(1+\alpha)^2+\ldots+(1+\alpha)^k)\\
    &\leq \beta(1+(1+\alpha)+(1+\alpha)^2+\ldots+(1+\alpha)^{2^m})\\
    &= \beta\frac{(1+\alpha)^{2^m+1}-1}{\alpha}\\
    &=\frac{a^2L(M+1)}{2^{2m+2}} \frac{(1+\frac{aL}{2^m})^{2^m+1}-1}{\frac{aL}{2^m}}\\
    &\leq\frac{a(M+1)\left(e^{aL}(1+\frac{aL}{2^m})-1\right)}{2^{m+2}}
\end{align*}
If $C=\frac{a(M+1)\left(e^{aL}(1+aL)-1\right)}{4}$, as the above expression is independent of $k$, it is deduced that, for all $d\in\D$ such that $e_d^m$ is well defined,
$$e_d^m\leq \frac{C}{2^m}.$$
From this it is deduced that the sequence $(x_d^m)_{m=m_d}^\infty$ is Cauchy: if $n_0\leq m<n$,
$$|x_d^m-x_d^n|\leq e_d^m+\ldots+ e_d^{n-1}\leq C\left(\frac{1}{2^m}+\ldots+\frac1{2^{n-1}}\right)<\frac{C}{2^{n_0-1}}$$
and therefore $x_d^m\to x_d$ for a certain $x_d\in \mathbb{R}$. In particular $|x_d^m-x_d|\leq \frac{C}{2^{m-1}}\,\forall d\in\D:m\geq m_d$ then the convergence is uniform in $d$. $\Box$

\begin{defi}
    For each $d\in\D$, we call $x_d=\lim_{m\to\infty}x_d^m$.
\end{defi}

\begin{lemma}
    For all $c,d\in\D$, $|x_c-x_d|\leq M|c-d|$.
\end{lemma}
\textbf{Proof.} It follows from Lemma 1.
\begin{lemma}
    If $(d_n)_{n\in\N}\subset\D$ is convergent, then $(x_{d_n})_{n\in\N}$ is convergent.
\end{lemma}
\textbf{Proof.} $(d_n)$ is Cauchy and from Lemma 2 it follows that $(x_{d_n})$ is Cauchy then convergent.
\begin{lemma}
    If $(c_n)_{n\in\N},(d_n)_{n\in\N}\subset\D$ are sequences convergent to the same point, then $(x_{c_n})_{n\in\N}$ and $(x_{d_n})_{n\in\N}$ are convergent to the same point.
\end{lemma}
\textbf{Proof.} By Lemma 3, $(x_{c_n})_{n\in\N}$ and $(x_{d_n})_{n\in\N}$ are convergent. Also, as $|x_{c_n}-x_{d_n}|\leq M|c_n-d_n|\to 0$, it follows that $|x_{c_n}-x_{d_n}|\to 0$.

\begin{defi}
    For each $t\in[0,a]$, we define $x(t)=\lim_{n\to \infty}x_{d_n}$ where $(d_n)_{n\in\N}\subset\D$ is any sequence convergent to $t$. Observe that $x(d)=x_d$ for all $d\in\D$.
\end{defi}
With Definition 4 we have obtained a function $x:[0,a]\to\mathbb{R}$. We will now see that it is a solution of the IVP.
\begin{lemma}
    For all $t,s\in[0,a]$, $|x(t)-x(s)|\leq M|t-s|$.
\end{lemma}
\textbf{Proof.} It follows from Lemma 2.
\begin{coro}
    For all $t\in[0,a]$, $(t,x(t))\in U$.
\end{coro}
\begin{prop}
If $(c_n)_{n\in\N},(d_n)_{n\in\N}\subset\D$ are sequences converging both to $t\in[0,a]$ such that $c_n<d_n\, \forall n\in\N$, then
$$\lim_{n\to\infty}\frac{x_{c_n}-x_{d_n}}{c_n-d_n}=f(t,x(t)).$$
\end{prop}
\textbf{Proof.}
\begin{equation}\left|\frac{x_{c_n}-x_{d_n}}{c_n-d_n}-f(t,x(t))\right|=\lim_{m\to\infty}\left|\frac{x_{c_n}^m-x_{d_n}^m}{c_n-d_n}-f(t,x(t))\right|.\end{equation}
Now, for each $m$ (sufficiently large) we can write $d_n=c_n+\frac{ak_n^m}{2^m}$. Then
$$x_{d_n}^m=x_{c_n}^m+\frac{a}{2^m}\left(f(c_n,x_{c_n}^m)+f(c_n+\frac{a}{2^m},x_{c_n+\frac{a}{2^m}}^m)+\ldots+f(d_n-\frac{a}{2^m},x_{d_n-\frac{a}{2^m}}^m)\right).$$
So,
\begin{align*}
&\left|\frac{x_{c_n}^m-x_{d_n}^m}{c_n-d_n}-f(t,x(t))\right|=\\
&=\left|\frac{\frac{a}{2^m}\left(f(c_n,x_{c_n}^m)+f(c_n+\frac{a}{2^m},x_{c_n+\frac{a}{2^m}}^m)+\ldots+f(d_n-\frac{a}{2^m},x_{d_n-\frac{}{2^m}}^m)\right)}{\frac{a}{2^m}k_n^m}-f(t,x(t))\right|\\
&=\left|\frac{\big(f(c_n,x_{c_n}^m)-f(t,x(t))\big)+\ldots+\big(f(d_n-\frac{a}{2^m},x_{d_n-\frac{a}{2^m}}^m)-f(t,x(t))\big)}{k_n^m}\right|\\
&\leq\frac{\left|f(c_n,x_{c_n}^m)-f(c_n,x_{c_n})\right|+\ldots+\left|f(d_n-\frac{a}{2^m},x_{d_n-\frac{a}{2^m}}^m)-f(d_n-\frac{a}{2^m},x_{d_n-\frac{a}{2^m}})\right|}{k_n^m}\\
&+\frac{\left|f(c_n,x_{c_n})-f(t,x(t))\right|+\ldots+\left|f(d_n-\frac{a}{2^m},x_{d_n-\frac{a}{2^m}})-f(t,x(t))\right|}{k_n^m}\\
&\leq L\frac{\left|x_{c_n}^m-x_{c_n}\right|+\ldots+\left|x_{d_n-\frac{a}{2^m}}^m-x_{d_n-\frac{a}{2^m}}\right|}{k_n^m}\qquad\text{(Lipschitz property)}\\
&+L(M+1)\frac{\left|c_n-t\right|+\ldots+\left|d_n-\frac{a}{2^m}-t\right|}{k_n^m}\qquad\text{(Lischitz property followed by Lemma 5).}
\end{align*}
Now, as $c_n<d_n$, it is easily checked that for all $0\leq k\leq k_n^m$, $\left|c_n+\frac{ak}{2^m}-t\right|\leq \max(|c_n-t|,|d_n-t|)=:\mu_n$. Also, as the convergence of $(x^m_d)$ is uniform in $d$, for all $\epsilon>0$ there exists $m_\epsilon$ so that if $m\geq m_\epsilon$ then $\left|x_{c_n}^m-x_{c_n}\right|,\ldots,\left|x_{d_n-\frac{a}{2^m}}^m-x_{d_n-\frac{a}{2^m}}\right|<\epsilon$. Combining this fact with the previous one, it is found that if $m\geq m_\epsilon$:
$$\left|\frac{x_{c_n}^m-x_{d_n}^m}{c_n-d_n}-f(t,x(t))\right|\leq L\frac{\epsilon+\ldots+\epsilon}{k_n^m}+L(M+1)\frac{\mu_n+\ldots+\mu_n}{k_n^m}.$$
As the summands of the numerators have exactly $k_n^m$ terms, and returning to (4), it is concluded that
$$\left|\frac{x_{c_n}-x_{d_n}}{c_n-d_n}-f(t,x(t))\right|=\lim_{m\to\infty}\left|\frac{x_{c_n}^m-x_{d_n}^m}{c_n-d_n}-f(t,x(t))\right|\leq L\epsilon+L(M+1)\mu_n.$$
As this occurs for all $\epsilon>0$, we obtain
$$\left|\frac{x_{c_n}-x_{d_n}}{c_n-d_n}-f(t,x(t))\right|\leq L(M+1)\mu_n.$$
As $\mu_n\to 0$, we finally conclude that
$$\lim_{n\to\infty}\left|\frac{x_{c_n}-x_{d_n}}{c_n-d_n}-f(t,x(t))\right|=0$$
as we wanted to see. $\Box$
\begin{prop}
    The function $x$ is differentiable and, for all $t\in[0,a]$, $x'(t)=f(t,x(t))$.
\end{prop}
\textbf{Proof.} We will first take $t\in ]0,1[$. Let $(t_n)_{n\in\N}$ be a sequence convergent to $t$. We first assume that $t<t_n\,\forall n\in\N$. Due to the density of $\D$, we can find sequences $(c_n)_{n\in\N},(d_n)_{n\in\N}\subset\D$ such that
\begin{itemize}
    \item $t<c_n<d_n<t_n$
    \item $\lim_{n\to\infty}\frac{d_n-c_n}{t_n-t}=1$.
\end{itemize}
In particular, we have $\lim_{n\to\infty}\frac{t_n-d_n}{t_n-t}=\lim_{n\to\infty}\frac{c_n-t}{t_n-t}=0$ and, as a consequence of Proposition 2, $$\lim_{n\to\infty}\frac{x_{d_n}-x_{c_n}}{t_n-t}=f(t,x(t)).$$

Then
\begin{align*}
    &\left|\frac{x(t_n)-x(t)}{t_n-t}-f(t,x(t))\right|\\
    &\leq \left|\frac{x(t_n)-x(d_n)}{t_n-t}\right|+\left|\frac{x(c_n)-x(t)}{t_n-t}\right|+\left|\frac{x(d_n)-x(c_n)}{t_n-t}-f(t,x(t))\right|\\
    &\leq M \left|\frac{t_n-d_n}{t_n-t}\right|+M\left|\frac{c_n-t}{t_n-t}\right|+\left|\frac{x(d_n)-x(c_n)}{t_n-t}-f(t,x(t))\right|\to0
\end{align*}
by the sequential characterization of the limit, $\lim_{s\to t^+}\frac{x(s)-x(t)}{s-t}=f(t,x(t))$. The limit from the left is calculated with an analogous method and gives $f(t,x(t))$. If $t=0$ or $t=1$ only the limit from the right and left respectively would be taken, according to the definition of derivative at those points. $\Box$

\begin{lemma}
    If $u:[0,b]\to\mathbb{R}$ is differentiable, $u(0)=0$ and for all $t\in[0,b]$ then $|u'(t)|\leq K|u(t)|$ for a certain real constant $K>0$, then $u(t)=0\,\forall t\in [0,b]$.
\end{lemma}
\textbf{Proof.} $u$ is continuous on a compact set then bounded, let's say $|u|\leq A>0$. Let $S=\{t\in[0,b]:u(t)=0\}$. $S$ is not empty. Let $s\in S$ and let $t_0\in [0,b]$ with $s<t_0<s+\frac1{2K}$. Applying the mean value theorem successively, $\{t_n\}_{n\geq 1}$ are obtained with $t_0>t_1>t_2>t_3>\ldots>t_n>\ldots>s$ so that for all $n\geq 0$
$$u(t_n)=u(t_n)-u(s)=u'(t_{n+1})(t_n-s)\implies |u(t_n)|\leq |u'(t_{n+1})|\frac{1}{2K}\leq \frac{|u(t_{n+1})|}{2}.$$
But then for all $n\geq 0$, $|u(t_0)|\leq \frac{|u(t_n)|}{2^n}\leq \frac{A}{2^n}$ from which it follows that $u(t_0)=0$.

This means that if $s\in S$ then $[s,s+\frac1{2K}[\cap [0,b]\subset S$. It is immediate that $S=[0,b]$ as we wanted to see. $\Box$

\begin{theor}
    Under the conditions mentioned on the first page about $f$, the IVP admits a unique solution with domain $[0,a]$.
\end{theor}
\textbf{Proof.} The existence is already proven. Let's now suppose that there exist $x_1,x_2:[0,a]\to\mathbb{R}$ such that $x_i(0)=x_0$ and, for all $t\in[0,a]$, $ (t,x_i(t))\in D$ and $x_i'(t)=f(t,x_i(t))$ for $i=1,2$. The first thing is to see that $(t,x_i(t))\in U\,\forall t\in [0,a]$, $i=1,2$.

Let $t=\min\{\tau\in[0,a[:|x_i(\tau)-x_0|\geq Ma\}>0$, if BWOC this set is non-empty. As $x_i$ is continuous, by the intermediate value theorem, there exists $s\in]0,t[$ such that $|x_i(s)-x_0|=M\frac{a+t}{2}>Mt$. By the mean value theorem there exists $0<w<s$ with
$$|f(w,x_i(w))| =|x_i'(w)|=\frac{|x_i(s)-x_0|}{s}>\frac{Mt}{t}=M$$
which implies that $(w,x_i(w))\notin U\implies |x_i(w)-x_0|>Ma$ contradicting the minimality of $t$.

With this we now know that for all $t\in [0,a]$ we have $(t,x_i(t))\in U$. Now, let $u=x_1-x_2$. Observe that $u$ is differentiable and $u(0)=0$. Also, for all $t\in [0,a]$
$$|u'(t)|=|x_1'(t)-x_2'(t)|=|f(t,x_1(t))-f(t,x_2(t))|\leq L|x_1(t)-x_2(t)|=L|u(t)|.$$
We have been able to apply the Lipschitz property since $(t,x_i(t))\in U$, $i=1,2$. By Lemma 6, $u=0$ and $x_1=x_2$. $\Box$

\end{document}